\documentclass[12pt, twoside, leqno]{article}
\usepackage{amsmath,amsthm}
\usepackage{amssymb,latexsym}
\usepackage{enumerate}
\usepackage[noadjust]{cite}

\newtheorem{thm}{Theorem}[section]

\newtheorem{cor}[thm]{Corollary}
\newtheorem{lem}[thm]{Lemma}

\theoremstyle{definition}

\newtheorem{rem}[thm]{Remark}

\numberwithin{equation}{section}

\usepackage{mathrsfs}
\usepackage{CJK}
\usepackage{fancyhdr}
\usepackage{amsmath}

\usepackage{enumerate}

\newcommand{\GL}{\operatorname{GL}}
\newcommand{\Mat}{\operatorname{Mat}}
\newcommand{\rk}{\operatorname{rank}}
\newcommand{\trdeg}{\operatorname{trdeg}}
\newcommand{\N}{\mathbb{N}}
\newcommand{\ideal}{\mathfrak{i}}

\begin{document}

\baselineskip=17pt

\title{The classification of some polynomial maps with nilpotent Jacobians}
\author{ Dan Yan \footnote{ The author is supported by the National Natural
Science Foundation of China (Grant No.11601146), the Natural Science
Foundation of Hunan Province (Grant No.2016JJ3085) and the Construct
Program of the Key Discipline in Hunan Province.}\\
MOE-LCSM,\\ School of Mathematics and Statistics,\\
 Hunan Normal University, Changsha 410081, China \\
\emph{E-mail:} yan-dan-hi@163.com \\
\and
Michiel de Bondt \\
Department of Mathematics, Radboud University \\
Nijmegen, The Netherlands \\
\emph{E-mail:} M.deBondt@math.ru.nl}
\date{}

\maketitle

\renewcommand{\thefootnote}{}

\renewcommand{\thefootnote}{\arabic{footnote}}
\setcounter{footnote}{0}

\begin{abstract} In the paper, we first classify all
polynomial maps $H$ of the following form:
$H=\big(H_1(x_1,x_2,\ldots,x_n),H_2(x_1,x_2),H_3(x_1,x_2),\ldots,H_n(x_1,x_2)\big)$
with $JH$ nilpotent. After that, we generalize the structure of $H$ to
$H=\big(H_1(x_1,x_2,\ldots,x_n),H_2(x_1,x_2),H_3(x_1,x_2,H_1),\ldots,H_n(x_1,x_2,H_1)\big)$.
\end{abstract}
 \vskip 2.5mm

\section{Introduction}

Throughout this paper, we will write $K$ for any field with
characteristic zero and $K[x]=K[x_1,x_2,\ldots,x_n]$ for the
polynomial algebra over $K$ with $n$ indeterminates. Let
$F=(F_1,F_2,\ldots,F_n):K^n\rightarrow K^n$ be a polynomial map,
that is, $F_i\in K[x]$ for all $1\leq i\leq n$. Let
$JF=(\frac{\partial F_i}{\partial x_j})_{n\times n}$ be the Jacobian
matrix of $F$. For $f\in K[x]$, we abbreviate $\frac{\partial
f}{\partial x_i}$ as $f_{x_i}$, and define $\deg_{x_i} f$ as the
highest degree of variable $x_i$ in $f$. $P_n(i,j)$ denotes the $n
\times n$ elementary permutation matrix which interchanges
coordinates $i$ and $j$.

The Jacobian Conjecture (JC) raised by Keller in 1939 in \cite{1}
states that a polynomial map $F:K^n\rightarrow K^n$ is invertible if
the Jacobian determinant $\det JF$ is a nonzero constant. This
conjecture has been attacked by many people from various research
fields, but it is still open, for all $n\geq 2$. Only the case $n=1$
is obvious. For more information about the wonderful 70-year
history, see \cite{2}, \cite{3}, and the references therein.

In 1980, Wang \cite{4} showed that the JC holds for all polynomial
maps of degree 2 in all dimensions. A powerful result is the
reduction to degree 3, due to Bass, Connell and Wright \cite{2} in
1982 and Yagzhev \cite{5} in 1980, which asserts that the JC is true
if it holds for all polynomial maps $x+H$, where $H$ is homogeneous
of degree 3. Thus, many authors studied these maps and which led to
pose the following problem.

{\em (Homogeneous) Dependence Problem.} Let $H=(H_1,\ldots,H_n)\in
K[x]$ be a (homogeneous) polynomial map of degree $d$ such that $JH$
is nilpotent and $H(0)=0$. Are $H_1,\ldots,H_n$ linearly dependent
over $K$?

The answer to the above problem is affirmative if $\rk JH\leq 1$
(\cite{2}). In particular, this implies that the Dependence Problem
has an affirmative answer in the case $n=2$. The second author and Van den
Essen give an affirmative answer to the above problem in the case that $H$
is homogeneous and $n=3$ (\cite{8}), which was extended by the second
author to the case $n=4$, under the assumption that $\rk JH\leq 2$
(\cite[Theorem 4.1.2]{HKM}).

With restrictions on the degree of $H$, more positive results are known.
For cubic homogeneous $H$, the case $n = 4$ has been solved affirmatively
by Hubbers in \cite{7}, using techniques of \cite{6}. For cubic homogeneous
$H$ with $\rk JH = 2$, the Dependence Problem has an affirmative answer for every
$n$, because the missing case $n \ge 5$ follows from \cite[Theorem 4.3.1]{HKM}.
For cubic $H$, the case $n = 3$ has been solved affirmatively as well,
see e.g.\@ \cite[Corollary 4.6.6]{HKM}).

For quadratic $H$, the Dependence Problem has an affirmative answer if
$\rk \allowbreak JH \le 2$ (see \cite{B2} or \cite[Theorem 3.4]{12}), in
particular if $n \le 3$. For quadratic homogeneous $H$,
the Dependence Problem has an affirmative answer in the case $n \le 5$,
and several authors contributed to that result. See \cite[Appendix A]{HKM}
and \cite{XS5} for the case $n = 5$. Recently, the second author generalized
the condition $n \le 5$ to $\rk JH\leq 4$ (\cite{B4}).

The first counterexamples to the Dependence Problem were found by
Van den Essen (\cite{9}, \cite[Theorem 7.1.7 (ii)]{3}). He
constructs counterexamples for all $n \ge 3$. In another paper
(\cite{E}), he constructs a quadratic counterexample for $n = 4$, which can be
generalized to arbitrary even degree (see \cite[Example 8.4.4]{3} for degree $4$).
Furthermore, the quadratic counterexample can be used to make a counterexample
of any degree $d \ge 2$ in larger dimensions:
$$
\big( x_2, x_1^2 - x_4, x_2^2, 2 x_1 x_2 - x_3, x_4^d, x_5^d, \ldots, x_{n-1}^d \big)
$$
We did not find a cubic counterexample for $n = 4$ in the literature. But
$$
\big( {-x_2^{d-1}} - (x_1 + x_2 x_3)x_3 + x_4, x_1 + x_2 x_3, x_2^{d-2}, (x_1 + x_2 x_3)x_2^{d-2} \big)
$$
is a counterexample of degree $d$ for every $d \ge 3$.

The second author was the
first who found homogeneous counterexamples (\cite{10}). He
constructed homogeneous counterexamples of $6$
for $n = 5$, homogeneous counterexamples of degree $4$ and $5$
for all $n \ge 6$, and cubic homogeneous counterexamples for all $n \ge 10$.
Homogeneous counterexamples of larger degrees can be made as well,
except for $n = 5$ and odd degrees. A cubic homogeneous counterexample
for $n = 9$ can be found in \cite{SFGZ}, see also \cite[Section 4.2]{HKM}.
We do not know any homogeneous counterexample of odd degree for $n = 5$.

In \cite{18}, Chamberland and Van den Essen classified
all polynomial maps of the form
$$H=\big(u(x_1,x_2),v(x_1,x_2,x_3),h(u(x_1,x_2),v(x_1,x_2,x_3))\big)$$
with $JH$ nilpotent. The first author and Tang \cite{13} classified
all polynomial maps of the form
$H=\big(u(x_1,x_2),v(x_1,x_2,x_3),h(x_1,x_2,x_3)\big)$ with some
conditions. In \cite{CE}, Casta\~{n}eda and Van den Essen classify all polynomial maps
of the form
$$H=\big(u(x_1,x_2),u_2(x_1,x_2,x_3),u_3(x_1,x_2,x_4),\ldots,u_{n-1}(x_1,x_2,x_n),
h(x_1,x_2)\big)$$ with $JH$ nilpotent. For more results about classification of
Keller maps, see e.g.\@ \cite{HKM,S3,B2,B3,B4,8,11,12,7,PC,14,XS5}. But the classification
result of \cite{7} (see also \cite[Theorem 7.1.2]{3}) is overly complicated,
see \cite[Theorem 4.6.5]{HKM} for the case where $K$ is algebraically closed,
and \cite{B3} for the general case.

A polynomial map of the form
$(x_1,\ldots,x_{i-1},x_i+P,x_{i+1},\ldots,x_n)$ is \emph{elementary} if
$P\in K[x_1,\ldots,x_{i-1},x_{i+1},\ldots,x_n]$. A polynomial map is
called \emph{tame} if it is a finite composition of invertible linear maps
and elementary maps.

In Section 2--3 of this paper, we will classify all polynomial maps
of the form
$H=\big(H_1(x_1,x_2),H_2(x_1,x_2,\ldots,x_n\big),H_3(x_1,x_2),\ldots,\allowbreak
H_n(x_1,x_2)\big)$ with $JH$ nilpotent, and show that the
corresponding Keller maps $x+H$ are tame. In fact, we will show that
$x+tH$ is tame over $K[t]$, where $t$ is a new variable. And in
Section 4--5, we will extend these results to the case where
$$
H=\big(H_1(x_1,x_2,\ldots,x_n),H_2(x_1,x_2),H_3(x_1,x_2,H_1),\ldots,H_n(x_1,x_2,H_1)\big).
$$

\section{The case where the components are linearly independent}

In this section, we classify all polynomial maps of the form
$$
H=\big(H_1(x_1,x_2,\ldots,\allowbreak x_n),H_2(x_1,x_2),H_3(x_1,x_2),\ldots, H_n(x_1,x_2)\big)
$$
such that $JH$ is nilpotent and the components of $H$ are linearly
independent over $K$. We additionally make the free assumption that
$H(0)=0$.

\begin{lem} \label{lm2.1}
Let $q \in K[x_1,x_2]$, such that $q_{x_1} \mid q_{x_2}^k$ for some
$k \in \N$, and $q_{x_1}$ does not have a factor in $K[x_2]
\setminus K$. Then there exists a polynomial $b \in K[x_2]$, such
that $q \in K[x_1 + b(x_2)]$.
\end{lem}

\begin{proof}
If $q_{x_2} = 0$, then $q \in K[x_1]$ and we can take $b = 0$.
So assume that $q_{x_2} \ne 0$.
%
We first show that
\begin{equation} \label{qx1qx2}
q_{x_1} \mid q_{x_2}.
\end{equation}
Let $p \in K[x_1,x_2]$ be an irreducible divisor of $q_{x_1}$.
Then $p \notin K[x_2]$ and $p \mid q_{x_2}$. In particular,
$(q_{x_1})_{x_2} = (q_{x_2})_{x_1} \ne 0$. Write $v_p(f)$ for
the multiplicity of $p$ as a divisor of a polynomial $f$. Then
$$
v_p(q_{x_1}) \le v_p\big((q_{x_1})_{x_2}\big) + 1 = v_p\big((q_{x_2})_{x_1}\big) + 1.
$$
Since $p \notin K[x_2]$ and $p \mid q_{x_2}$, we infer that
$$
v_p\big((q_{x_2})_{x_1}\big) + 1 = v_p(q_{x_2}),
$$
so $v_p(q_{x_1}) \le v_p(q_{x_2})$. This holds for every
irreducible $p \mid q_{x_1}$, which yields equation \eqref{qx1qx2}.

Since $\deg_{x_1} (q_{x_2}) \le \deg_{x_1} q = \deg_{x_1} (q_{x_1}) + 1$,
we infer from equation \eqref{qx1qx2} that
\begin{equation} \label{qx2qx1}
q_{x_2} = \big(a(x_2) x_1 + b'(x_2)\big) q_{x_1}
\end{equation}
for univariate polynomials $a, b$.

Let $u = x_1^i x_2^j$ be the term of highest degree with respect to $x_2$
among the terms of highest degree with respect to $x_1$ of $q$
(so $i = \deg_{x_1} q$). If $a \ne 0$, then the right-hand side of equation
\eqref{qx2qx1} has degree $i$ with respect to $x_1$, just like
$u$, so $j \ge 1$ and
$$
\frac{u}{x_2} = x_2^{\deg a} x_1 \frac{u}{x_1},
$$
which is impossible. So $a = 0$ and $q_{x_2} = b'(x_2) q_{x_1}$.
Consequently, $D q = 0$, where
$D := \partial_{x_2}-b'(x_2)\partial_{x_1}$.
Notice that $D \big(x_1 + b(x_2)\big) = 0$ and
$D x_2 = 1$. Since $q \in K[x_1,x_2] = K[x_1 + b(x_2),x_2]$,
we infer that $q \in K[x_1 + b(x_2)]$.
\end{proof}

We will not use the following generalization of Lemma \ref{lm2.1},
but it is used in \cite{CE}: \cite[Lemma 2]{CE} is a slightly different
formulation of Corollary \ref{cor2.1} below. \cite[Lemma 2]{CE} is proved
by reducing to the case where $K$ is algebraically closed. This special
case of \cite[Lemma 2]{CE} appeared in an earlier version of this article
by the first author, and \cite{CE} refers to that. Below, we give a direct
proof.

\begin{cor} \label{cor2.1}
Let $q \in K[x_1,x_2]$, such that $q_{x_1} \mid w\,q_{x_2}^k$ for some
$w \in K[q]$ and $k \in \N$, and $q_{x_1}$ does not have a factor in
$K[x_2] \setminus K$. Then there exists a polynomial $b \in K[x_2]$,
such that $q \in K[x_1 + b(x_2)]$.
\end{cor}

\begin{proof}
Let $Q \in K[q]$ and $p \in K[x_2] \setminus K$. If $p \mid Q_{x_1} = Q_q\,q_{x_1}$,
then $p \mid Q_q$ because $\gcd (p, q_{x_1}) = 1$, and $p \mid (Q_q)_{x_1}$ because
$p_{x_1} = 0$. By induction on $\deg_q Q$, we obtain that $p \nmid Q_{x_1}$.
So $Q_{x_1}$ does not have a factor in $K[x_2] \setminus K$.

Choose $Q \in K[q]$, such that $Q_q = w$. Then
$$
Q_{x_1} = Q_q\,q_{x_1} = w\,q_{x_1} \mid w^2 q_{x_2}^k \mid w^{k+2} q_{x_2}^{k+2} = Q_{x_2}^{k+2}
$$
On account of lemma \ref{lm2.1}, $Q \in K[x_1 + b(x_2)]$ for some $b \in K[x_2]$.
Hence for $D := \partial_{x_2}-b'(x_2)\partial_{x_1}$, we have $Dq = DQ / Q_q = 0$,
and $q \in K[x_1 + b(x_2)]$ follows in a similar manner as in the end of
the proof of lemma \ref{lm2.1}.
\end{proof}

\begin{lem} \label{lmeq}
Let $H$ be a polynomial map over $K$ of the form
$$
\big(H_1(x_1,x_2,x_3\ldots,x_n),H_2(x_1,x_2),H_3(x_1,x_2),\ldots,H_n(x_1,x_2)\big).
$$
If $JH$ is nilpotent, then
\begin{align}
(H_1)_{x_1}+(H_2)_{x_2} &= 0, \label{eq2.3} \\
\big((H_2)_{x_2}\big)^2 + (H_1)_{x_2}(H_2)_{x_1} + (H_1)_{x_3}(H_3)_{x_1} +
\cdots + (H_1)_{x_n}(H_n)_{x_1} &= 0, \label{eq2.4} \\
\begin{split}
(H_1)_{x_3}\big((H_2)_{x_1}(H_3)_{x_2}-(H_2)_{x_2}(H_3)_{x_1}\big) +
(H_1)_{x_4}\big((H_2)_{x_1}(H_4)_{x_2}-{} \\
(H_2)_{x_2}(H_4)_{x_1}\big) + \cdots + (H_1)_{x_n}\big((H_2)_{x_1}(H_n)_{x_2}-
(H_2)_{x_2}(H_n)_{x_1}\big) &= 0
\end{split} \label{eq2.5}
\end{align}
and $(H_2)_{x_1} \mid \big((H_2)_{x_2}\big)^3$.
\end{lem}

\begin{proof}
Suppose that $JH$ is nilpotent.
Adding equation \eqref{eq2.4} $(H_2)_{x_2}$ times to equation \eqref{eq2.5}
yields the last claim. So it remains to prove equations \eqref{eq2.3},
\eqref{eq2.4} and \eqref{eq2.5}.

Equation \eqref{eq2.3} follows from the fact that the trace of $JH$ is zero.
Since the sum of the principal minor determinants of size $2$ of $JH$ is zero
as well, we deduce that
$$
-(H_1)_{x_1}(H_2)_{x_2}+(H_1)_{x_2}(H_2)_{x_1}+(H_1)_{x_3}(H_3)_{x_1}+
\cdots+(H_1)_{x_n}(H_n)_{x_1} = 0.
$$
Adding equation \eqref{eq2.3} $(H_2)_{x_2}$ times to it yields equation
\eqref{eq2.4}. Equation \eqref{eq2.5} follows from the fact that the sum of
the principal minor determinants of size $3$ of $JH$ is zero.
\end{proof}

\begin{lem} \label{lm2.2}
Let $H$ be a polynomial map over $K$ of the form
$$
\big(H_1(x_1,x_2,x_3\ldots,x_n),H_2(x_1,x_2),H_3(x_1,x_2),\ldots,H_n(x_1,x_2)\big).
$$
Assume that $H(0) = 0$, and that the components of $H$ are linearly
independent over $K$. If $JH$ is nilpotent, then $\deg \bar{H}_1=1$,
where $\bar{H}_1$ is the leading homogeneous part with respect to
$x_3, x_4,\ldots, x_n$ of $H_1$.
\end{lem}

\begin{proof}
Suppose that $JH$ is nilpotent. We first prove that
$H_1 \notin K[x_1,x_2]$. Suppose that $H_1 \in K[x_1,x_2]$.
Then the nilpotency of $JH$ comes down to the
nilpotency of the leading principal minor $P$
of size $2$ of $JH$.

Since $P \in \Mat_2(K[x_1,x_2])$,
we deduce that $P$ is a nilpotent jacobian matrix itself.
On account of e.g.\@ Theorem 7.1.7 i) in \cite{3}, the rows of
$P$ are dependent over $K$. This contradicts the fact that the
components of $H$ are independent over $K$, so
$H_1 \notin K[x_1,x_2]$.

Let $\bar{H}_1$ be the leading homogeneous part with respect to
$x_3, x_4, \ldots, x_n$ of $H_1$. Then $\bar{H}_1$ has positive
degree, say $d$, with respect to $x_3, x_4, \ldots, x_n$. Our aim is
to prove that $\bar{H}_1$ is a linear combination of $x_3, x_4,
\ldots, x_n$. Notice that, as far as it is nonzero,
$(\bar{H}_1)_{x_i}$ has degree $d$ with respect to $x_3, x_4,
\ldots, x_n$ if $i \le 2$, and degree $d-1$ with respect to $x_3,
x_4, \ldots, x_n$ otherwise.

If we focus on the leading homogeneous part with respect to
$x_3, x_4, \ldots, x_n$ of equation \eqref{eq2.3}, we deduce that
$(\bar{H}_1)_{x_1} = 0$. If $(H_2)_{x_1} = 0$, then
$(H_2)_{x_2}$ would be an eigenvalue of $JH$, which is impossible
because $H_2 \notin K$. So $(H_2)_{x_1} \ne 0$. If we focus on the
leading homogeneous part with respect to $x_3, x_4, \ldots, x_n$
of equation \eqref{eq2.4}, we deduce that $(\bar{H}_1)_{x_2} = 0$.
So $\bar{H}_1 \in K[x_3,x_4,\ldots,x_n]$.
Assume without loss of generality that
$$
(\bar{H}_1)_{x_3},(\bar{H}_1)_{x_4},\ldots,(\bar{H}_1)_{x_k}
$$
are linearly independent over $K$, and $(\bar{H}_1)_{x_{k+1}} =
(\bar{H}_1)_{x_{k+2}} = \cdots = (\bar{H}_1)_{x_n} = 0$.

Then $(\bar{H}_1)_{x_3},(\bar{H}_1)_{x_4},\ldots,
(\bar{H}_1)_{x_k}$ are linearly independent over $K(x_1,x_2)$
as well. So if we focus on the leading homogeneous part with
respect to $x_3, x_4, \ldots, x_n$ of equation \eqref{eq2.5}, we infer
that
$$
(H_2)_{x_1} (H_i)_{x_2} -(H_2)_{x_2} (H_i)_{x_1} = 0
$$
for each $i \in \{3,4,\ldots,k\}$. Consequently,
$H_i$ is algebraically dependent over $K$ on $H_2$ for
each $i \in \{3,4,\ldots,k\}$, and there exists an
$f \in K[x_1,x_2]$, such that $H_i \in K[f]$
for each $i \in \{2,3,4,\ldots,k\}$. As $(H_2)_{x_1} \ne 0$,
we have $f \in K[x_1,x_2] \setminus K[x_2]$.

Now assume that $d > 1$. Let $H^{*}_1$ be the homogeneous part
of degree $d-1 > 0$ with respect to
$x_3, x_4, \ldots, x_n$ of $H_1$.
If we focus on the leading homogeneous part with respect to
$x_3, x_4, \ldots, x_n$ of equation \eqref{eq2.3} again, we deduce that
$(H^{*}_1)_{x_1} = 0$. Notice that $(H_i)_{x_1}/(H_2)_{x_1} \in
K(f)$ for each $i \ge 3$.
So if we focus on the leading homogeneous part with respect to
$x_3, x_4, \ldots, x_n$ of equation \eqref{eq2.4} again, we infer that
$$
(H^{*}_1)_{x_2} \in K(f)(\bar{H}_1)_{x_3} + K(f)(\bar{H}_1)_{x_4} +
\cdots + K(f)(\bar{H}_1)_{x_n} \subseteq K(f)[x_3,x_4,\ldots,x_n].
$$
From $(H^{*}_1)_{x_1} = 0$ above, we deduce that
$(H^{*}_1)_{x_2} \in K[x_2,x_3,x_4,\ldots,x_n]$ and
that $\big((H^{*}_1)_{x_2}\big)_{x_1} = 0$. But $f_{x_1} \ne 0$, so
every coefficient of $(H^{*}_1)_{x_2} \in K(f)[x_3,x_4,\ldots,x_n]$
is contained in $K$, i.e.\@
$(H^{*}_1)_{x_2} \in K[x_3,x_4,\ldots,x_n]$. So if
we focus on the leading homogeneous part with
respect to $x_3, x_4, \ldots, x_n$ of equation \eqref{eq2.4} for the third
time, we infer that
$$
(H_2)_{x_1},(H_3)_{x_1},\ldots,(H_k)_{x_1}
$$
are linearly dependent over $K(x_3,x_4,\ldots,x_n)$, and hence over $K$.
Since the rank of the submatrix of rows $2,3,\ldots,k$ of $JH$ is $1$,
the rows of this submatrix are (linearly) dependent over $K$ along with
the entries of its first column. This contradicts the fact that the
components of $H$ are linearly independent over $K$, so $\deg \bar{H}_1 = d = 1$.
\end{proof}

\begin{thm} \label{th2.3}
Let $H$ be a polynomial map over $K$ of the form
$$
\big(H_1(x_1,x_2,x_3\ldots,x_n),H_2(x_1,x_2),H_3(x_1,x_2),\ldots,H_n(x_1,x_2)\big).
$$
Assume that $H(0) = 0$, and that the components of $H$ are linearly
independent over $K$. If $JH$ is nilpotent,
then there are
$b_1, b_2 \in K$, such that $H_2 \in K[x_1 + b_2 x_2^2 + b_1 x_2]$,
and $\sigma_2, \sigma_3, \ldots, \sigma_{n} \in K$ such that
\begin{align*}
H_1 + (2 b_2 x_2 + b_1)H_2 &=
\sigma_2 x_2 + \sigma_3 x_3 + \cdots + \sigma_n x_n, \\
b_2 H_2^2 &= \sigma_2 H_2 + \sigma_3 H_3 + \cdots + \sigma_n H_n.
\end{align*}
Furthermore, one of $\sigma_3, \sigma_4, \ldots, \sigma_n$
is nonzero, as well as $b_2$.
\end{thm}

\begin{proof}
It follows from Lemma \ref{lm2.2} that $\deg\bar{H}_1=1$, where $\bar{H}_1$ is the
leading homogeneous part with respect to $x_3, x_4,\ldots, x_n$ of $H_1$.
Write
$$
\bar{H}_1 = \sigma_3 x_3 + \sigma_4 x_4 + \cdots + \sigma_n x_n
$$
with $\sigma_3, \sigma_4, \ldots, \sigma_n \in K$, and define
$h := \sigma_3 H_3 + \sigma_4 H_4 + \cdots + \sigma_n H_n$.
Then equation \eqref{eq2.5} comes down to
\begin{align}
0 = \sum_{i=3}^n \sigma_i \big((H_2)_{x_1}(H_i)_{x_2} - (H_2)_{x_2}
(H_i)_{x_1}\big) &= (H_2)_{x_1} h_{x_2} - (H_2)_{x_2} h_{x_1}.
\label{eq2.6}
\end{align}
So $h$ is algebraically dependent over $K$ on $H_2$, and there
exists an
$f \in K[x_1,x_2]$, such that $H_2, h \in K[f]$. Say that
$H_2 = g(f)$, where $g$ is an univariate polynomial over $K$.
From equations \eqref{eq2.4} and \eqref{eq2.3}, we infer that
\begin{align}
\begin{split}
0 &= \big((H_2)_{x_2}\big)^2 + (H_1)_{x_2} (H_2)_{x_1} +
\sum_{i=3}^n \sigma_i (H_i)_{x_1} \\
&= -(H_1)_{x_1}(H_2)_{x_2} + (H_1)_{x_2} (H_2)_{x_1} + h_{x_1}.
\end{split} \label{eq2.7}
\end{align}

We define the leading term of a polynomial $p \in K[x_1,x_2]$ as
the term of highest degree with respect to $x_2$ among the terms
of highest degree with respect to $x_1$ of $p$.
Now let $u_1$, $u_2$, and $u_3$ be the leading terms of
$H_1 - \bar{H}_1$, $H_2$, and $h$ respectively.
Write $u_2 = x_1^i x_2^j$. We distinguish two cases:
\begin{itemize}

\item $j \ge 1$.

Then it follows from equation \eqref{eq2.3} that $u_1 = x_1^{i+1}x_2^{j-1}$.
Therefore,
$$
(u_1)_{x_1}(u_2)_{x_2} - (u_1)_{x_2} (u_2)_{x_1} = \big((i+1)j -
(j-1)i\big) x_1^{2i}x_2^{2j-2} \ne 0.
$$
If we focus on the coefficient of $x_1^{2i}x_2^{2j-2}$ of
equation \eqref{eq2.7}, we see that the leading term of $h_{x_1}$ is
$x_1^{2i}x_2^{2j-2}$. Hence $u_3 = x_1^{2i+1}x_2^{2j-2}$.
Just as for $u_1$ and $u_2$ before,
$$
(u_2)_{x_1} (u_3)_{x_2} - (u_2)_{x_2} (u_3)_{x_1} = \big(i(2j-2) -
j(2i+1)\big) x_1^{3i}x_2^{3j-3} \ne 0.
$$
This contradicts the fact that the coefficient of
$x_1^{3i}x_2^{3j-3}$ of equation \eqref{eq2.6} is zero.

\item $j = 0$.

From Lemma \ref{lm2.1} and the last claim of Lemma \ref{lmeq},
it follows that $H_2 \in K[x_1+b(x_2)]$ for some univariate
polynomial $b$ over $K$, such that $b(0) = 0$.
Hence we can take $f = x_1 + b(x_2)$, so $\deg g = i$.
Furthermore, $g(0) = 0$ because $H(0) = 0$.
From equation \eqref{eq2.3}, we deduce that
$H_1$ is of the form
$$
-b'(x_2)g(f) + c(x_2) + \sigma_3 x_3 + \sigma_4 x_4 + \cdots + \sigma_n x_n
$$
for an univariate polynomial $c$ over $K$. As $H(0) = 0$,
we see that $c(0) = 0$.

Notice that $h_{x_1} \in K[f]$ because $h \in K[f]$. From
equation \eqref{eq2.7}, it follows that
$(H_1)_{x_1} (H_2)_{x_2} - (H_2)_{x_1} (H_1)_{x_2} \in K[f]$.
Since
\begin{align*}
(H_1)_{x_1} &= -b'(x_2) g'(f) &
(H_1)_{x_2} &= -b''(x_2) g(f) - b'(x_2)^2 g'(f) + c'(x_2) \\
(H_2)_{x_1} &= g'(f) & (H_2)_{x_2} &= b'(x_2) g'(f),
\end{align*}
we infer that $-b''(x_2) g(f) + c'(x_2) \in K(f)$
(the terms with $b'(x_2)$ cancel out).
Hence $-b''(x_2) g(x_1) + c'(x_2) \in K(x_1)$ and
$-b'''(x_2) g(x_1) + c''(x_2) = 0$. So
$b'''(x_2) = c''(x_2) = 0$, and we can write
$b = b_2 x_2^2 + b_1 x_2$ and $c = \sigma_2 x_2$.
It follows from equation \eqref{eq2.7} that
$$
h_{x_1} = (H_1)_{x_1} (H_2)_{x_2} - (H_2)_{x_1} (H_1)_{x_2} = 2 b_2
g'(f) g(f) - \sigma_2 g'(f).
$$
Consequently, $h = b_2 g(f)^2 - \sigma_2 g(f)$.
Now the equalities in Theorem \ref{th2.3} can be obtained from
the formulas for $H_2$, $H_1$, and $h$. \qedhere

\end{itemize}
\end{proof}

\section{The case where the components may be linearly dependent}

In this section, we remove the assumption of the previous section
that the components of $H$ are linearly independent over $K$. More
generally, we classify all maps $H$ with $JH$ nilpotent, such that
$H_i \in K[x_1,x_2]$ for all $i$ with one exception, and show that
$x+H$ is tame.

\begin{thm} \label{th2.4}
Let $H$ be a polynomial map over $K$ of the form
$$
\big(H_1(x_1,x_2,\ldots,x_n), H_2(x_1,x_2),H_3(x_1,x_2),
\ldots,H_n(x_1,x_2)\big),$$
such that $H(0) = 0$. If $JH$ is
nilpotent, then there exists a $T \in \GL_n(K)$, such that
$\tilde{H} = T^{-1}H(Tx)$ has the same form as $H$ itself, and one
of the following statements holds:
\begin{enumerate}[\upshape (i)]

\item $\tilde{H}_2 = 0$, and there exists a $k \in \{2,3,\ldots,n\}$
and an $A \in \GL_{n-2}(K[x_2])$, such that $\tilde{H}_1 \in
K[x_2,A_{k-1} \tilde{x}, A_k \tilde{x},\ldots, A_{n-2} \tilde{x}]$, and
$$
\left( \begin{smallmatrix}
\tilde{H}_3(x_1,x_2) - \tilde{H}_3(0,x_2) \\
\tilde{H}_4(x_1,x_2) - \tilde{H}_4(0,x_2) \\
\vdots \\
\tilde{H}_n(x_1,x_2) - \tilde{H}_n(0,x_2)
\end{smallmatrix} \right)
$$
is a $K[x_1,x_2]$-linear combination of the first $k-2$ columns of
$A^{-1}$, where $A_l$ denotes the $l$-th row of $A$, and
$\tilde{x} = (x_3,x_4,\ldots,x_n)$ as a column vector.

\item $\tilde{H}_2 \in K[x_1] \setminus K$ and there exists a
$k \in \{2,3,\ldots,n\}$, such that
$\tilde{H}_1 \in K[x_{k+1},x_{k+2},\ldots,x_n]$, and
$\tilde{H}_{k+1} = \tilde{H}_{k+2} = \cdots = \tilde{H}_{n} = 0$.

\item $\tilde{H}_2 \in K[x_1+x_2^2]$, $\tilde{H}_3 = \tilde{H}_2^2$,
and there exists a $k \in \{3,4,\ldots,n\}$, such that
$$
\tilde{H}_1 + 2x_2 \tilde{H}_2 - x_3 \in K[x_{k+1},x_{k+2},\ldots,x_n]
$$
and $\tilde{H}_{k+1} = \tilde{H}_{k+2} = \cdots = \tilde{H}_{n} = 0$.

\end{enumerate}
\end{thm}

\begin{proof}
Take $T \in \GL_n(K)$, such that $\tilde{H} = T^{-1}H(Tx)$
has the same form as $H$ itself.
We distinguish two cases.
\begin{itemize}

\item \emph{$H_1$ and $H_2$ are linearly dependent over $K$.}

Then we can choose $T$, such that $\tilde{H}_2 = 0$. Hence
equation \eqref{eq2.5} is satisfied for $\tilde{H}$.
Equation \eqref{eq2.3} for $\tilde{H}$ tells us that the upper left
corner of $J\tilde{H}$ is zero along with the rest of the
diagonal of $J\tilde{H}$. So equation \eqref{eq2.4} for $\tilde{H}$
comes down to that the first row of $J\tilde{H}$ is orthogonal
to the first column of $J\tilde{H}$.
$x_2$ is the only variable which is in both.

Let $d = \deg_{x_1} (\tilde{H}_3,\tilde{H}_4,\ldots,\tilde{H}_n)$,
and let $m \in K[x_1,x_2]^{n-2}$ be the vector which consists of
entries $3,4,\ldots,n$ of the first column of $J\tilde{H}$.
Then we can write
\begin{equation} \label{m}
m = M \cdot \left( \begin{smallmatrix}
1 \\ x_1 \\ x_1^2 \\ \vdots \\ x_1^{d-1}
\end{smallmatrix} \right),
\end{equation}
where $M \in \Mat_{n-2,d}(K[x_2])$. Now antidifferentiating with
respect to $x_1$ yields
\begin{equation} \label{mH}
\left( \begin{smallmatrix}
\tilde{H}_3(x_1,x_2) - \tilde{H}_3(0,x_2) \\
\tilde{H}_4(x_1,x_2) - \tilde{H}_4(0,x_2) \\
\vdots \\
\tilde{H}_n(x_1,x_2) - \tilde{H}_n(0,x_2)
\end{smallmatrix} \right) = M \cdot \left( \begin{smallmatrix}
x_1 \\ \frac12 x_1^2 \\ \frac13 x_1^3 \\ \vdots \\ \frac1d x_1^d
\end{smallmatrix} \right),
\end{equation}
Let $k = \rk M + 2$. Since
$K[x_2]$ is an Euclidean domain, we can reduce $M$ with row
operations, so there exists an $A \in \GL_{n-2}(K[x_2])$ such that
only the first $k-2$ rows of $AM$ are nonzero. Since $M = A^{-1}
AM$, we deduce that the left hand side of \eqref{mH} is a
$K[x_1,x_2]$-linear combination of the first $k-2$ columns of
$A^{-1}$.

As $\tilde{H}_1\big(x_1,x_2,(A^{-1})_1 \tilde{x},(A^{-1})_2 \tilde{x},
\ldots, (A^{-1})_{n-2} \tilde{x}\big) \in K[x_2,x_3,x_4,\ldots,x_n]$,
we deduce that
\begin{equation} \label{eq2.9}
\tilde{H}_1 \in K[x_2,A_1 \tilde{x}, A_2 \tilde{x}, \ldots, A_{n-2}
\tilde{x}].
\end{equation}
Since $\rk M = k-2$, we infer that the column space of $M$
over $K(x_2)$ is equal to the space over $K(x_2)$ generated by
the first $k-2$ columns of $A^{-1}$.
From \eqref{m} and the fact that the first row of $J\tilde{H}$ is orthogonal to
the first column of $J\tilde{H}$, we deduce that
$$
\big((A^{-1})_{1j} \partial_3 + (A^{-1})_{2j} \partial_4 +
\cdots + (A^{-1})_{(n-2)j} \partial_n\big) \tilde{H}_1 = 0
$$
for all $j \le k-2$. So by \eqref{eq2.9},
$\tilde{H}_1 \in K[x_2,A_{k-1} \tilde{x}, A_k \tilde{x},\ldots,
A_{n-2} \tilde{x}]$, and $\tilde{H}$ is as in (i) above.

\item \emph{$H_1$ and $H_2$ are linearly independent over $K$.}

Since $H_1$ and $H_2$ are linearly independent over $K$,
we can choose $T$ such that
$\tilde{H}_1, \tilde{H}_2, \tilde{H}_3, \ldots, \tilde{H}_{k}$
are linearly independent over $K$, and
$\tilde{H}_{k+1} = \tilde{H}_{k+2} = \cdots = \tilde{H}_{n} = 0$.

Suppose first that $\tilde{H}_1 \in K[x_1,x_2,x_{k+1},x_{k+2},
\ldots,x_n]$. Then the nilpotency of $J \tilde{H}$ comes down to the
nilpotency of the leading principal minor
$$
P = \left( \begin{array}{cc} P_{11} & P_{12} \\
P_{21} & P_{22} \end{array} \right)
$$
of size $2$ of $J\tilde{H}$. If $P_{21} = 0$, then $P_{22}$
would be an eigenvalue of $P$, which is impossible.
So $P_{21} \ne 0$.

$P_{12}$ is the only entry of $P$
which is not necessarily contained in $K[x_1,x_2]$. Since $P$
is nilpotent and $P_{21} \ne 0$, we deduce that
$P \in \Mat_2(K[x_1,x_2])$.
Thus, $P$ is a nilpotent jacobian matrix itself.
On account of e.g.\@ Theorem 7.1.7 i) in \cite{3}, the rows of
$P$ are dependent over $K$.
Thus, $P_{21}$ and $P_{11}$ are linearly dependent over $K$.
Since $P_{21} \ne 0$, we can choose $T$ such that $P_{11} = 0$.
Since $P$ is nilpotent, we infer that
$P_{12} = P_{22} = 0$. Hence $\tilde{H}_1 \in
K[x_{k+1},x_{k+2},\ldots,x_n]$, $\tilde{H}_2 \in K[x_1]$,
and $\tilde{H}$ is as in (ii) above.

Suppose next that $\tilde{H}_1 \notin K[x_1,x_2,x_{k+1},x_{k+2},
\ldots,x_n]$. Then $\tilde{H}_1, \tilde{H}_2, \tilde{H}_3,
 \ldots, \allowbreak\tilde{H}_{k}$ are linearly independent
over $K(x_{k+1},x_{k+2},\ldots,x_n)$.
Furthermore, the nilpotency of $J\tilde{H}$ comes down to the
nilpotency of the leading principal minor of size $k$ of
$J\tilde{H}$, so we can apply Theorem \ref{th2.3}, with base field
$K(x_{k+1},x_{k+2},\ldots,x_n)$ instead of $K$, and
dimension $k$ instead of $n$.

So $\tilde{H}_2 = g(x_1+b_2 x_2^2 + b_1 x_2)$ for some univariate
polynomial $g$ over $K(x_{k+1},\allowbreak x_{k+2},\ldots,x_n)$, and
$b_2, b_1 \in K(x_{k+1},x_{k+2},\ldots,x_n)$, with $b_2 \ne 0$.
But $\tilde{H}_2 \in K[x_1,x_2]$, and by substituting $x_2 = 0$,
we see that $g$ is just an univariate polynomial over $K$, say
of degree $i$. If we differentiate $\tilde{H}_2$ $i-1$ times
with respect to $x_1$, we infer that $b_2, b_1 \in K$.

Furthermore, $\tilde{H}_2^2$ can be written as a
linear combination over
$K(x_{k+1},\allowbreak x_{k+2},\allowbreak \ldots,x_n)$ of
$\tilde{H}_2,\tilde{H}_3,\ldots,\tilde{H}_{k}$, and hence as
a linear combination over $K$ as well. Since the coefficient of
one of $\tilde{H}_3,\ldots,\tilde{H}_{k}$ is nonzero, we can
choose $T$ such that $\tilde{H}_3 = \tilde{H}_2^2$.

Notice that
$$
\tilde{H}_2(b_2 x_1 - b_1 x_2,x_2) \in K[b_2x_1 - b_1x_2 + b_1 x_2 +
b_2 x_2^2] = K[x_1 + x_2^2].
$$
Hence we can choose $T$ such that $\tilde{H}_2 \in K[x_1 + x_2^2]$.
This way, we get $b_2 = 1$ and $b_1 = 0$,

If in our application of Theorem \ref{th2.3} with base field
$K(x_{k+1},x_{k+2},\allowbreak\ldots,x_n)$ and $k$ instead of $n$,
if either one of $\sigma_2, \sigma_4, \sigma_5, \ldots, \sigma_k$ is
nonzero, or $\sigma_3 \ne 1$, then $\tilde{H}_2, \tilde{H}_3,
\ldots, \tilde{H}_{k}$ would be linearly dependent over
$K(x_{k+1},x_{k+2},\ldots, \allowbreak x_n)$. So $\sigma_2 =
\sigma_4 = \sigma_5 = \cdots = \sigma_k = 0$, $\sigma_3 = 1$, and
$\tilde{H}$ is as in (iii) above. \qedhere

\end{itemize}
\end{proof}

\begin{rem}
The interpretation of Theorem \ref{th2.4} in the case where $k = n$ is, that
$K[x_2,A_{k-1} \tilde{x}, A_k \tilde{x}, \ldots,A_{n-2} \tilde{x}] = K[x_2]$,
that $K[x_{k+1},x_{k+2},\ldots,x_n] = K$, and that
$\tilde{H}_{k+1} = \tilde{H}_{k+2} = \cdots = \tilde{H}_{n} = 0$ is void.
\end{rem}

\begin{cor} \label{cr2.5}
Let $H = (H_1,H_2,\ldots,H_n)$ be a polynomial map over $K$,
such that for some $i \le n$, $H_j \in K[x_1,x_2]$ for all
$j \ne i$. Suppose that $H(0) = 0$ and that $JH$ is nilpotent.
\begin{itemize}

\item If $i \le 2$, then there exists a $T \in \GL_n(K)$, such
that $T^{-1} H(Tx)$ is the same as in
Theorem \ref{th2.4}.

\item If $i \ge 3$, then there exists a $T \in \GL_n(K)$, such
that $T^{-1} H(Tx)$ is of the form
\begin{equation} \label{eq2.10}
\big(0,\tilde{H}_2(x_1),\tilde{H}_3(x_1,x_2),\ldots,
\tilde{H}_{n-1}(x_1,x_2),\tilde{H}_n(x_1,x_2,\ldots,x_{n-1})\big).
\end{equation}

\end{itemize}
Furthermore, $x + tH$ is tame over $K[t]$. In particular,
$x + \lambda H$ is invertible over $K$ for every $\lambda \in K$.
\end{cor}

\begin{proof} Let $T \in \GL_n(K)$ and
$\tilde{H} := T^{-1}H(Tx)$.
\begin{itemize}

\item Suppose that $i \le 2$. Then there exists an elementary
permutation matrix $S$, such that $S^{-1}H(Sx)$ is as $H$
in Theorem \ref{th2.4}. Hence Theorem \ref{th2.4} holds for $H$.

\item Suppose that $i \ge 3$. Then we can take an elementary
permutation matrix for $T$, such that $\tilde{H}_j \in K[x_1,x_2]$
for every $j < n$.
As a consequence, $(\tilde{H}_n)_{x_n}$ is an eigenvalue of
$J\tilde{H}$. So $\tilde{H}_n \in K[x_1,x_2,\ldots,x_{n-1}]$,
and the nilpotency of $JH$ comes down to that of the leading
principal minor $P$ of size $2$ of $JH$. Just as in the proof
of Theorem \ref{th2.4}, we can obtain that the lower left corner of $P$ is
the only nonzero entry of $P$. So $\tilde{H} = T^{-1}H(Tx)$ is
of the given form.

\end{itemize}
So it remains to prove that $x + tH$ is tame over $K[t]$.
It suffices to show that $x + t\tilde{H}$ is tame over $K[t]$.
Since $J\tilde{H}$ is lower triangular in \eqref{eq2.10}, the case
$i \le 2$ remains. We distinguish the three cases of Theorem \ref{th2.4}:
\begin{enumerate}[\upshape (i)]

\item Since $K[x_2]$ is a Euclidean domain, we can write
$A$ as a product of elementary matrices over $K[x_2]$.
Hence $G = (x_1,x_2,A\tilde{x})$ is tame. Now a straightforward
computation yields that
$$
G\big(G^{-1} + t \tilde{H}(G^{-1})\big) = x + t H^{*}
$$
with $H^{*} = (H^{*}_1, H^{*}_2, \ldots, H^{*}_n)$ such that
$H^{*}_1 \in K[x_2,x_{k+1},x_{k+2},\ldots,x_n]$, $H^{*}_2 = 0$,
$H^{*}_j \in K[x_1,x_2]$ for all $j \ge 3$, and
$H^{*}_j \in K[x_2]$ for all $j \ge k+1$.

Since $x_2$ is a component of $x + t H^{*}$, we can get rid of terms
in $K[x_2]$ of $\tilde{H}_j$ for any $j \ne 2$ by way of an
elementary invertible map over $K[t]$ from the left. If we do this
for all $j \ge k+1$, then we can get a polynomial map $x +
t\hat{H}$, such that $\hat{H}$ is of the same form as $\tilde{H}$ in
Theorem \ref{th2.4} (ii), except that $\hat{H}_2 = 0 \ne
\tilde{H}_2$ and $\hat{H}_1\in K[x_2,x_{k+1},x_{k+2},\ldots,x_n]$.
Since $x_2, x_{k+1},x_{k+2},\ldots,x_n$ are components of $x + t
H^{*}$, we can get rid of terms in
$K[x_2,x_{k+1},x_{k+2},\ldots,x_n]$ of $\hat{H}_1$ by way of an
elementary invertible map over $K[t]$ from the left. So we may
assume that $\hat{H}_1$ has no terms in
$K[x_2,x_{k+1},x_{k+2},\ldots,x_n]$. Hence $\hat{H}_1=0$ and
$J\hat{H}$ is lower triangular. So $x + t\hat{H}$ is tame over
$K[t]$. And so is $x + t \tilde{H}$.

\item Since $x_{k+1},x_{k+2},\ldots,x_n$ are components of
$x + t \tilde{H}$, we can get rid of terms in
$K[x_{k+1},x_{k+2},\ldots,x_n]$ of $\tilde{H}_1$ by way of an
elementary invertible map over $K[t]$ from the left.
So we may assume that
$\tilde{H}_1$ has no terms in $K[x_{k+1},x_{k+2},\ldots,\allowbreak
x_n]$. Hence $\tilde{H}_1 = 0$ and $J\tilde{H}$ is lower triangular.
So $x + t \tilde{H}$ is tame over $K[t]$.

\item Just as above, we may assume that $\tilde{H}_1$ has no terms
in $K[x_{k+1},x_{k+2},\ldots,\allowbreak x_n]$. Hence
\begin{align*}
x_1 + t \tilde{H}_1 &= x_1 + t (-2 x_2 \tilde{H}_2 + x_3) \\
&= x_1 + x_2^2 - x_2^2 + t (-2 x_2 \tilde{H}_2 - t\tilde{H}_2^2 + x_3 + t \tilde{H}_2^2)\\
&= (x_1 + x_2^2) - (x_2 + t \tilde{H}_2)^2 + t(x_3 + t \tilde{H}_3).
\end{align*}
So by applying an elementary map over $K[t]$ from the left,
and by applying the elementary map
$(x_1-x_2^2,x_2,x_3,\ldots,x_n)$ from the right,
we obtain a polynomial map in lower triangular form from
$x + t \tilde{H}$. Hence $x + t \tilde{H}$ is tame over $K[t]$.
\qedhere

\end{enumerate}
\end{proof}

\section{A generalization of the form of {\mathversion{bold}$H$}}

From now on, $H$ is a polynomial map over $K$ of the form
$$
\big(H_1(x_1,x_2,x_3\ldots,x_n),H_2(x_1,x_2),h_3(x_1,x_2,H_1),\ldots,h_n(x_1,x_2,H_1)\big)
$$
such that $JH$ is nilpotent. The following lemma is an analog of Lemma \ref{lmeq}

\begin{lem}
\begin{align*}
\begin{split}
(H_1)_{x_1}+(H_2)_{x_2} + (h_3)_{x_3}(x_1,x_2,H_1) (H_1)_{x_3} + {} \\
(h_4)_{x_3}(x_1,x_2,H_1) (H_1)_{x_4} + \cdots + (h_n)_{x_3}(x_1,x_2,H_1) (H_1)_{x_n} &= 0
\end{split} \tag{\ref{eq2.3}$'$} \label{eq2.3a} \\
\begin{split}
\big((H_2)_{x_2}\big)^2 + (H_1)_{x_2}(H_2)_{x_1} + (H_1)_{x_3}(h_3)_{x_1}(x_1,x_2,H_1) + {} \\
(H_1)_{x_4}(h_4)_{x_1}(x_1,x_2,H_1) + \cdots + (H_1)_{x_n}(h_n)_{x_1}(x_1,x_2,H_1) &= 0
\end{split} \tag{\ref{eq2.4}$'$} \label{eq2.4a} \\
\begin{split}
(H_1)_{x_3}\big((H_2)_{x_1}(h_3)_{x_2}(x_1,x_2,H_1)-
(H_2)_{x_2}(h_3)_{x_1}(x_1,x_2,H_1)\big) + \cdots + \\
(H_1)_{x_n}\big((H_2)_{x_1}(h_n)_{x_2}(x_1,x_2,H_1)-
(H_2)_{x_2}(h_n)_{x_1}(x_1,x_2,H_1)\big) &= 0
\end{split} \tag{\ref{eq2.5}$'$} \label{eq2.5a}
\end{align*}
and $(H_2)_{x_1} \mid \big((H_2)_{x_2}\big)^3$.
\end{lem}

\begin{proof}
We can derive equation \eqref{eq2.3a},
\eqref{eq2.4a} and \eqref{eq2.5a} in a similar manner as equation
\eqref{eq2.3}, \eqref{eq2.4} and \eqref{eq2.5} respectively, provided
we can find proper formulas for the principal minor determinants of size
$2$ and $3$. Furthermore $(H_2)_{x_1} \mid \big((H_2)_{x_2}\big)^3$ can
be obtained just as before.

One way to find proper formulas for the principal minor determinants of
size $2$ and $3$ is to use $(H_1)_{x_j} (H_1)_{x_i} = (H_1)_{x_i} (H_1)_{x_j}$
to cancel out terms. But it can be done as follows as well. The determinant
of a minor matrix with at least $2$ columns with index larger than $2$ is
zero, because the last $2$ columns are dependent. To compute the determinant
of a minor matrix with row $1$ and one row $i \ge 3$, we can clean the
$i$-th row by means of the first row, and compute the corresponding
minor determinant of the result:
$$
\begin{pmatrix}
(H_1)_{x_1} & (H_1)_{x_2} & (H_1)_{x_3} & \cdots & (H_1)_{x_n} \\
(H_2)_{x_1} & (H_2)_{x_2} & 0 & \cdots & 0 \\
(h_i)_{x_1}(x_1,x_2,H_1) & (h_i)_{x_2}(x_1,x_2,H_1) & 0 & \cdots & 0
\end{pmatrix}
$$
Other principal minor determinants can be computed directly.
\end{proof}

By substituting $t = H_1$, we see that equations \eqref{eq2.3b},
\eqref{eq2.4b} and \eqref{eq2.5b} below, in which $t$ is a new variable,
are stronger than equations \eqref{eq2.3a}, \eqref{eq2.4a} and
\eqref{eq2.5a} above, respectively.
\begin{align*}
\begin{split}
(H_1)_{x_1} + (H_2)_{x_2} + (h_3)_{x_3}(x_1,x_2,t) (H_1)_{x_3} + {} \\
(h_4)_{x_3}(x_1,x_2,t) (H_1)_{x_4} + \cdots + (h_n)_{x_3}(x_1,x_2,t)
(H_1)_{x_n} &= 0,
\end{split} \tag{\ref{eq2.3a}$'$} \label{eq2.3b} \\
\begin{split}
\big((H_2)_{x_2}\big)^2 + (H_1)_{x_2}(H_2)_{x_1} + (H_1)_{x_3}(h_3)_{x_1}(x_1,x_2,t) + {} \\
(H_1)_{x_4}(h_4)_{x_1}(x_1,x_2,t) + \cdots +
(H_1)_{x_n}(h_n)_{x_1}(x_1,x_2,t) &= 0,
\end{split} \tag{\ref{eq2.4a}$'$} \label{eq2.4b} \\
\begin{split}
(H_1)_{x_3}\big((H_2)_{x_1}(h_3)_{x_2}(x_1,x_2,t)-
(H_2)_{x_2}(h_3)_{x_1}(x_1,x_2,t)\big) + \cdots + \\
(H_1)_{x_n}\big((H_2)_{x_1}(h_n)_{x_2}(x_1,x_2,t)-
(H_2)_{x_2}(h_n)_{x_1}(x_1,x_2,t)\big) &= 0.
\end{split} \tag{\ref{eq2.5a}$'$} \label{eq2.5b}
\end{align*}
We will prove equations \eqref{eq2.3b}, \eqref{eq2.4b} and \eqref{eq2.5b}
in the next section.

Whereas \eqref{eq2.4b} and \eqref{eq2.5b} seem good analogs of
\eqref{eq2.4} and \eqref{eq2.5} respectively, \eqref{eq2.3b}
does not seem a good analog of \eqref{eq2.3}, which is
\begin{equation*}
(H_1)_{x_1} + (H_2)_{x_2} = 0. \tag{\ref{eq2.3}}
\end{equation*}
Actually, there does not seem to be a good analog of \eqref{eq2.3} other than
\eqref{eq2.3}. More explicitly, if we have \eqref{eq2.3} instead of \eqref{eq2.3b},
then we can replace the argument $H_1$ of $h_3, h_4, \ldots, h_n$ by $t$
already in $H$ itself, instead of \eqref{eq2.3a}, \eqref{eq2.4a}, and
\eqref{eq2.5a}.

Unfortunately, \eqref{eq2.3} does not hold in general. But it is valid under
certain conditions.

\begin{lem} \label{lmeq2.3}
\eqref{eq2.3} is satisfied if one of the following statements holds.
\begin{enumerate}[\upshape(i)]

\item $H_2 \in K$ and $h_i \in \ideal(x_1,x_3^2)$ for each $i \ge 3$.

\item $H_2 \notin K$ and $h_i \in \ideal(x_1,x_2,x_3^2)$ for each $i \ge 3$.

\end{enumerate}
Here, $\ideal(g_1,g_2,\ldots)$ is the ideal in $K[x]$ which is generated by
$g_1,g_2,\ldots$.
\end{lem}

\begin{proof}[Proof of {\upshape(i)}]
Suppose that $H_2 \in K$, but $(H_1)_{x_1} + (H_2)_{x_2} \ne 0$.
Then $(H_1)_{x_1} \ne 0$. We will derive a contradiction.
Since $H_2 \in K$, we infer that equations \eqref{eq2.3b},
\eqref{eq2.4b} and \eqref{eq2.5b} remain valid if we replace $H_1$ by $(H_1)_{x_i}$
for any $i \ge 3$. Hence we may assume that $(H_1)_{x_1} \in K[x_1,x_2] \setminus \{0\}$.
So
$$
(H_1)_{x_i} \in K[x_2,x_3,\ldots,x_n]
$$
for each $i \ge 3$.

Suppose that $h_i \in \ideal(x_1,x_3^2)$ for each $i \ge 3$. Then
$t^2 \mid h_i(0,x_2,t)$ for each $i \ge 3$, and
we can antidifferentiate the left-hand side of equation
\eqref{eq2.4b} with respect to $x_1$, to obtain
$$
t^2 \mid (H_1)_{x_3}h_3(x_1,x_2,t) + (H_1)_{x_4}h_4(x_1,x_2,t) +
\cdots + (H_1)_{x_n}h_n(x_1,x_2,t).
$$
If we take the coefficient of $t^1$ of the right-hand side, then we
obtain the coefficient of $t^0$ of the left-hand side of equation
\eqref{eq2.3b}, except its first summand $(H_1)_{x_1}$. This
contradicts $(H_1)_{x_1} \ne 0$.
\end{proof}

The proof of Lemma \ref{lmeq2.3} (ii) will be given in the next section.

\begin{thm}
Assume that the components of $(H,1)$ are linearly independent over $K$.
Assume in addition that Lemma \ref{lmeq2.3} (ii) is satisfied.
Then there are $b_1, b_2 \in K$, such that
$H_2 \in K[x_1 + b_2 x_2^2 + b_1 x_2]$, and
$\sigma_2, \sigma_3, \ldots, \sigma_{n} \in K$ such that
\begin{alignat*}{2}
\big(H_1 + (2 b_2 x_2 + b_1)H_2\big) &-
\big(\sigma_2 x_2 + \sigma_3 x_3 + \cdots + \sigma_n x_n\big) && \in K, \\
b_2 H_2^2 &- \big(\sigma_2 H_2 + \sigma_3 H_3 + \cdots + \sigma_n
H_n\big) && \in K.
\end{alignat*}
If $H(0) = 0$ in addition, then the conclusions of Theorem \ref{th2.3}
hold for $H$, i.e.\@ we can replace $\in K$ by $= 0$ above.
\end{thm}

\begin{proof}
We can follow the proof of Theorem \ref{th2.3}, using
equation \eqref{eq2.3} as it is, and using \eqref{eq2.4b} and
\eqref{eq2.5b} as analogs of \eqref{eq2.4} and \eqref{eq2.5}
respectively.
\end{proof}

For a matrix $T \in \GL_n(K)$, we associate $T$ with the polynomial map
$Tx$: the matrix product of $T$ with $x$ as a column vector. So
$T \mapsto Tx$ is the embedding of $T \in \GL_n(K)$ into the $n$-dimensional
polynomial automorphisms. For a matrix $S \in \GL_{n-1}(K[x_2])$, $S \mapsto Sx$
does not work, because matrix dimensions do not match. But we can use
the following embedding
$$
\left( \begin{matrix}
S_{11} & S_{12} & \cdots & S_{1(n-1)} \\
S_{21} & S_{22} & \cdots & S_{2(n-1)} \\
\vdots & \vdots & \ddots & \vdots \\
S_{(n-1)1} & S_{(n-1)2} & \cdots & S_{(n-1)(n-1)}
\end{matrix} \right) \mapsto
\left( \begin{matrix}
S_{11} & 0 & S_{12} & \cdots & S_{1(n-1)} \\
0 & 1 & 0 & \cdots & 0 \\
S_{21} & 0 & S_{22} & \cdots & S_{2(n-1)} \\
\vdots & \vdots & \vdots & \ddots & \vdots \\
S_{(n-1)1} & 0 & S_{(n-1)2} & \cdots & S_{(n-1)(n-1)}
\end{matrix} \right) x.
$$

\begin{thm}
Suppose that $H_1 \notin K[x_1,x_2]$.
\begin{enumerate}[\upshape(i)]

\item If $H_2 = 0$, then there exists an $S \in \GL_{n-1}(K[x_2])$,
such that for $\tilde{H} = S^{-1}\big(H(Sx)\big)$,
\begin{align*}
\tilde{H}_1 &\in K[x_2,x_{k+1},x_{k+2},\ldots,x_n], \\
\tilde{H}_2 &= 0,
\end{align*}
and $\tilde{H}_{k+1}, \tilde{H}_{k+2}, \ldots, \tilde{H}_n \in K[x_2]$.

If $H_2 \in K^{*}$, then we have to prepare $H$ by removing the constant
part of $H_2$.

\item If $H_2 \in K[x_1] \setminus K$, then there exists a
$T \in \GL_{n-1}(K)$, such that for $\tilde{H} = T^{-1}\big(H(Tx)\big)$,
\begin{align*}
\tilde{H}_1 &\in K[x_{k+1},x_{k+2},\ldots,x_n], \\
\tilde{H}_2 &\in K[x_1] \setminus K,
\end{align*}
and $\tilde{H}_{k+1}, \tilde{H}_{k+2}, \ldots, \tilde{H}_n \in K$.

\item If $H_2 \in K[x_1,x_2] \setminus K[x_1]$, then there exists a
$T \in \GL_{n}(K)$, such that for $\tilde{H} = T^{-1}\big(H(Tx)\big)$,
\begin{align*}
\tilde{H}_1 + 2 x_2 \tilde{H}_2 - x_3 &\in K[x_{k+1},x_{k+2},\ldots,x_n], \\
\tilde{H}_2 &\in K[x_1 + x_2^2] \setminus K, \\
\tilde{H}_3 - \tilde{H}_2^2 &\in K,
\end{align*}
and $\tilde{H}_{k+1}, \tilde{H}_{k+2}, \ldots, \tilde{H}_n \in K$.

\end{enumerate}
Furthermore, $x + tH$ is tame over $K[t]$. In particular, $x + \lambda H$ is
invertible over $K$ for every $\lambda \in K$.
\end{thm}

\begin{proof}
We only prove (i), because (ii) and (iii) follow in a similar manner.
Suppose that $H_2 \in K$. By way of a conjugation process over
$\GL_{n-1}(K[x_2])$, we can get rid of terms $x_2^j x_3$ of $h_3, h_4,
\ldots, h_n$. But the condition of Lemma \ref{lmeq2.3} (i) instructs
that terms of the form $x_2^j$ of $h_3, h_4, \ldots, h_n$ should be canceled
as well. It is however obvious that terms of the form $x_2^j$ of
$h_3, h_4, \ldots, h_n$ do not have any influence on the left-hand sides
of equations \eqref{eq2.3b}, \eqref{eq2.4b} and \eqref{eq2.5b}. Hence
we can use Lemma \ref{lmeq2.3} (i) after all to justify the assumption that
\eqref{eq2.3} holds.

Now we can use the proof of Theorem \ref{th2.4}, to obtain that
there exists a $T \in \GL_{n-1}(K)$, such that for
$\tilde{H} = T^{-1}\big(H(Tx)\big)$, the following holds:
$\tilde{H}_2 = 0$, and there exists a $k \in \{2,3,\ldots,n\}$
and an $A \in \GL_{n-2}(K[x_2])$, such that $\tilde{H}_1 \in
K[x_2,A_{k-1} \tilde{x}, A_k \tilde{x},\ldots, A_{n-2} \tilde{x}]$, and
$$
\left( \begin{smallmatrix}
\tilde{h}_3(x_1,x_2,x_3) - \tilde{h}_3(0,x_2,x_3) \\
\tilde{h}_4(x_1,x_2,x_3) - \tilde{h}_4(0,x_2,x_3) \\
\vdots \\
\tilde{h}_n(x_1,x_2,x_3) - \tilde{h}_n(0,x_2,x_3)
\end{smallmatrix} \right)
$$
is a $K[x_1,x_2,x_3]$-linear combination of the first $k-2$ columns of
$A^{-1}$, where $A_l$ denotes the $l$-th row of $A$,
$\tilde{x} = (x_3,x_4,\ldots,x_n)$ as a column vector, and
$\tilde{h}_i$ is obtained by $\tilde{H}_i = \tilde{h}_i(x_1,x_2,\tilde{H}_1)$
for each $i \ge 3$.

Just as $h_3, h_4, \ldots, h_n$, $\tilde{h}_3, \tilde{h}_4, \ldots, \tilde{h}_n$
do not have terms of the form $x_2^j x_3$. To obtain $S \in \GL_{n-1}(K[x_2])$
as claimed, we can now follow the first step of (i) in the proof of Corollary
\ref{cr2.5}. The tameness of $x + t H$ over $K[t]$ follows in a similar manner
as before.
\end{proof}

\section{\mathversion{bold}Proof of equations (\ref{eq2.3b}), (\ref{eq2.4b}) and
(\ref{eq2.5b}), and Lemma \ref{lmeq2.3} (ii)\mathversion{normal}}

Let $w(f)$ be the weighted degree of a polynomial $f \in K[x,t]$,
such that $w(x_1) = w(x_2) = w(t) = 0$ and $w(x_i) \ge 1$ for all $i \ge
3$ (where $w(x_i)$ are real numbers). Let $\Phi(f)$ be the sum of the
monomials $\tau$ of $f \in K[x]$, for which $w(\tau) \ge w(H_1)$.
Then $\Phi(H_1)$ is the $w$-leading part of $H_1$, and
$\Phi\big((H_1)_{x_i}\big) = 0$ for all $i \ge 3$. Furthermore,
$\Phi(f)_{x_1} = \Phi(f_{x_1})$ and $\Phi(f)_{x_2} = \Phi(f_{x_2})$.

\begin{lem} \label{lmeqb}
\begin{enumerate}[\upshape(i)]

\item If $\Phi(H_1)_{x_1} = 0$, then equation \eqref{eq2.3b} holds.

\item If either $H_2 \in K$ or $\Phi(H_1)_{x_2} = 0$, then equation
\eqref{eq2.4b} holds.

\item \eqref{eq2.5b} holds.

\end{enumerate}
\end{lem}

\begin{proof}
The conditions of (i), (ii), and (iii) ensure that $\Phi\big((H_1)_{x_i}\big) = 0$
for all $i$ for which $(H_1)_{x_i}$ appears in equations \eqref{eq2.3b}, \eqref{eq2.4b},
and \eqref{eq2.5b} respectively. So if we take for $f$ the left-hand side of any equation
of equations \eqref{eq2.3b}, \eqref{eq2.4b}, and \eqref{eq2.5b}, then $\Phi(f) = 0$.

Let $w^{*}(f) = w(f) + w(H_1) \deg_t f$. Then one can infer from $\Phi(f) = 0$ that
$$
w^{*}(f) = w(f|_{t=H_1}).
$$
But $f|_{t=H_1}$ is the left-hand side of equations \eqref{eq2.3a}, \eqref{eq2.4a}, and
\eqref{eq2.5a}, which is equal to zero. Hence $w^{*}(f) = w(0) = -\infty$ and $f = 0$.
\end{proof}

For a polynomial $h \in K[x_1,x_2,x_3]$, we define
$h^{*}$ as the quotient polynomial of $h$ and $x_3$, i.e.\@
$$
h^{*}(x_1,x_2,x_3) := \frac{h(x_1,x_2,x_3) - h(x_1,x_2,0)}{x_3}.
$$
Furthermore,
\begin{equation}
h_{x_1}^{*} := (h_{x_1})^{*} = (h^{*})_{x_1},
 \quad h_{x_2}^{*} :=
(h_{x_2})^{*} = (h^{*})_{x_2} \quad \mbox{and} \quad h_{x_3}^{*} :=
(h_{x_3})^{*}.
\end{equation}
Notice that $h_{x_3}^{*} = (h^{*})_{x_3}$, if and only if
$\deg_{x_3} h \le 1$, if and only if either side is zero.

\begin{lem} \label{PhiH1}
\begin{enumerate}[\upshape(i)]

\item
$\Phi(H_1)$ is contained in $K[x_2,x_3,x_4,\ldots,x_n]$, and equation
\eqref{eq2.3b} holds.

\item
If $(H_2)_{x_1} \ne 0$, then $\Phi(H_1)$ is the product of a polynomial
in $K[x_2]$ and a polynomial in $K[x_3,x_4,\ldots,x_n]$.

\end{enumerate}
\end{lem}

\begin{proof}
\begin{enumerate}[\upshape(i)]

\item
From equation \eqref{eq2.3a}, it follows that
\begin{gather*}
\Phi\big((h_3)_{x_3}(x_1,x_2,H_1) (H_1)_{x_3} + \cdots
+ (h_n)_{x_3}(x_1,x_2,H_1) (H_1)_{x_n}\big) \\
= \Phi\big({-(H_1)_{x_1}}\big) = -\Phi(H_1)_{x_1}.
\end{gather*}
Subtracting monomials $\tau$ for which $w(\tau) < w(H_1)$ will
not affect that, so if we define
$$
\alpha_3 := (h_3)_{x_3}^{*}(x_1,x_2,H_1) (H_1)_{x_3} + \cdots +
(h_n)_{x_3}^{*}(x_1,x_2,H_1) (H_1)_{x_n},
$$
then $\Phi(H_1 \alpha_3) = -\Phi(H_1)_{x_1}$ as well.
Consequently, $w(H_1 \alpha_3) \le w(H_1)$, so $w(\alpha_3) \le 0$,
and
\begin{equation}
\Phi(H_1) \alpha_3 = \Phi(H_1 \alpha_3) = -\Phi(H_1)_{x_1}.
\end{equation}
If $\alpha_3 \ne 0$, then comparing degrees yields a contradiction,
so $\alpha_3 = 0$ and $\Phi(H_1)_{x_1} = 0$. Hence
$\Phi(H_1) \in K[x_2,x_3,x_4,\ldots,x_n]$, and equation \eqref{eq2.3b}
follows from Lemma \ref{lmeqb} (i).

\item
From equation \eqref{eq2.4a}, it follows that
\begin{gather*}
\Phi\big((H_1)_{x_3}(h_3)_{x_1}(x_1,x_2,H_1) + \cdots
+ (H_1)_{x_n}(h_n)_{x_1}(x_1,x_2,H_1)\big) \\
= \Phi\big({-(H_1)_{x_2}(H_2)_{x_1}}\big) =
-\Phi(H_1)_{x_2}(H_2)_{x_1}.
\end{gather*}
Subtracting monomials $\tau$ for which $w(\tau) < w(H_1)$ will
not affect that, so if we define
\begin{equation} \label{alpha}
\alpha_1 := (H_1)_{x_3}(h_3)_{x_1}^{*}(x_1,x_2,H_1) + \cdots +
(H_1)_{x_n}(h_n)_{x_1}^{*}(x_1,x_2,H_1),
\end{equation}
then $\Phi(H_1 \alpha_1) = -\Phi(H_1)_{x_2}(H_2)_{x_1}$ as well.
Consequently, $w(H_1 \alpha_1) \le w(H_1)$, so $w(\alpha_1) \le 0$,
and
\begin{equation} \label{alpha1}
\Phi(H_1) \alpha_1 = \Phi(H_1 \alpha_1) =
-\Phi(H_1)_{x_2}(H_2)_{x_1}.
\end{equation}
Let $p \in K[x_2,x_3,x_4,\ldots,x_n] \setminus K[x_3,x_4,\ldots,x_n]$
be a prime factor of $\Phi(H_1)$. Then $\Phi(H_1)_{x_2} \ne 0$. Write
$v_p(f)$ for the multiplicity of $p$ as a factor of $f$.
Then
$$
v_p\big(\Phi(H_1)\big) = v_p\big(\Phi(H_1)_{x_2}\big) + 1.
$$
From equation \eqref{alpha1}, we infer that $p \mid (H_2)_{x_1} \in
K[x_1,x_2] \setminus \{0\}$. Consequently, $p \in K[x_2]$.
From (i), we deduce that $\Phi(H_1)$ is as claimed. \qedhere

\end{enumerate}
\end{proof}

If $(H_2)_{x_1} = 0$, then equation \eqref{eq2.4b} follows from
Lemma \ref{lmeqb} (ii). On account of that and the condition in Lemma
\ref{lmeq2.3} (ii), we assume from now on that
$$
(H_2)_{x_1} \ne 0.
$$
We will prove equation \eqref{eq2.4b} and Lemma \ref{lmeq2.3} (ii) under
a special condition. After that, we show that the special condition is not
a real restriction.

\begin{lem} \label{histar}
Suppose that there exists a $k$ such that $H_1 \in K[x_1,x_2,\ldots,x_k]$,
and such that
$$
\Phi(H_1)_{x_3}, \Phi(H_1)_{x_4}, \ldots, \Phi(H_1)_{x_k}
$$
are linearly independent over $K$. Then
$\Phi(H_1)_{x_3}, \Phi(H_1)_{x_4}, \ldots, \Phi(H_1)_{x_k}$ are
linearly independent over $K(x_1,x_2)$, and the following hold.
\begin{enumerate}[\upshape(i)]

\item There exists an $f \in K[x_1,x_2] \setminus K[x_2]$, such that $H_2 \in K[f]$
and $h_i^{*} \in K[f]$ for all $i \in \{3,4,\ldots,k\}$.

\item $\Phi(H_1) \in K[x_3,x_4,\ldots,x_k]$, and equation \eqref{eq2.4b}
holds.

\item Lemma \ref{lmeq2.3} (ii) holds, i.e.\@ if
$h_i \in \ideal(x_1,x_2,x_3^2)$ for each $i \ge 3$, then
$(H_1)_{x_1} + (H_2)_{x_2} = 0$.

\end{enumerate}
\end{lem}

\begin{proof}
From Lemma \ref{PhiH1} (ii), we infer that
$$
\Phi(H_1)_{x_3}, \Phi(H_1)_{x_4}, \ldots, \Phi(H_1)_{x_k}
$$
are linearly independent over $K(x_1,x_2)$. Hence no top term cancelation can occur
in a linear combination of $(H_1)_{x_3}, (H_1)_{x_4}, \ldots, (H_1)_{x_k}$ over
$K(x_1,x_2)$, so
$$
(H_1)_{x_3}, (H_1)_{x_4}, \ldots, (H_1)_{x_k}
$$
are linearly independent over $K(x_1,x_2)$ as well.
\begin{enumerate}[\upshape(i)]

\item If we apply the above on the coefficients of $t^1, t^2, t^3, \ldots$
of the left-hand side of equation \eqref{eq2.3b}, we infer that
$\deg_{x_3} h_i \le 1$ for each $i \in \{3,4,\ldots,k\}$.
Hence $h_i^{*} \in K[x_1,x_2]$ for each $i \in \{3,4,\ldots,k\}$.

Furthermore, if we apply the above on the coefficient of $t^1$ of the left-hand
side of equation \eqref{eq2.5b}, we infer that the coefficient of $t^1$ of
$$
(H_2)_{x_1} (h_i)_{x_2}(x_1,x_2,t) - (H_2)_{x_2} (h_i)_{x_1}(x_1,x_2,t)
$$
is zero for all $i \in \{3,4,\ldots,k\}$. Consequently,
$\trdeg_K K(H_2,h_3^{*},\allowbreak h_4^{*},\ldots,h_k^{*}) = 1$.
Hence an $f \in K[x_1,x_2]$ as claimed exists, and $f \notin K[x_2]$
because $(H_2)_{x_1} \ne 0$.

\item Again, we use \eqref{alpha} to define $\alpha_1$.
From (i), it follows that $\alpha_1$ is a linear combination over
$K(x_1,x_2)$ of $(H_1)_{x_3}, \allowbreak (H_1)_{x_4}, \ldots, (H_1)_{x_k}$.
Just as above, no top term cancelation can occur, so
$$
w\big((H_1)_{x_i}(h_i)_{x_1}^{*}(x_1,x_2,H_1)\big) =
w\big((H_1)_{x_i}(h_i)_{x_1}^{*}\big) \le w(\alpha_1) \le 0
$$
for each $i \in \{3,4,\ldots,k\}$.
From Lemma \ref{lmeqb} (ii), we deduce that there exists a $\beta \in K[x_2]$,
such that $\Phi(H_1)/\beta \in K[x_3,x_4,\ldots,x_n]$. Now suppose that
$(h_i)_{x_1}^{*} \ne 0$ for some $i \in \{3,4,\ldots,k\}$.
Then $w\big((H_1)_{x_i}\big) = 0$. So $w(x_i) = w(H_1)$ and
\begin{equation} \label{fracconst}
\frac{(H_1)_{x_i}}{\beta} = \frac{\Phi(H_1)_{x_i}}{\beta} \in K.
\end{equation}
If we combine this with (i), we deduce that for each $i \in \{3,4,\ldots,k\}$,
$$
\frac{(H_1)_{x_i}(h_i)_{x_1}^{*}(x_1,x_2,H_1)}{\beta f_{x_1}} \in K[f],
\qquad \mbox{so} \qquad
\frac{\alpha_1}{\beta f_{x_1}} \in K[f].
$$
Write $H_2 = g(f)$. Then $(H_2)_{x_1} = f_{x_1} g'(f)$.
From equation \eqref{alpha1}, we infer that
$$
\frac{\alpha_1}{\beta f_{x_1}}
= -\frac{\Phi(H_1)_{x_2} (H_2)_{x_1}}{\Phi(H_1) \beta f_{x_1}}
= -\frac{\beta_{x_2} g'(f)}{\beta^2},
\qquad \mbox{so} \qquad
\frac{\beta_{x_2}}{\beta^2} \in K(f).
$$
Suppose that $\beta \notin K$. By definition of $\beta$, $f$ and $x_2$ are
algebraically dependent over $K$. This contradicts $f \notin K[x_2]$.
So $\beta \in K$, and the claims follow from the definition
of $\beta$ and Lemma \ref{lmeqb} (ii).

\item From equation \eqref{fracconst}, we deduce that
$$
(h_i)_{x_1}^{*} \ne 0 \Longrightarrow \Phi(H_1)_{x_i} \in K
$$
for each $i \in \{3,4,\ldots,k\}$. Hence $(h_i)_{x_1}^{*} \ne 0$ for
at most one $i \in \{3,4,\ldots,k\}$. So $\alpha_1$ has at most one
nonzero summand. From equation \eqref{alpha1} and $\Phi(H_1)_{x_2} = 0$,
we deduce that $\alpha_1 = 0$, so every summand of $\alpha_1$ is zero,
and $(h_i)_{x_1}^{*} = 0$ for all $i \in \{3,4,\ldots,k\}$.
As $h_i^{*} \in K[f]$ and $f \in K[x_1,x_2] \setminus K[x_1]$ for each
$i \in \{3,4,\ldots,k\}$, we deduce that $h_i^{*} \in K$ for every
$i \in \{3,4,\ldots,k\}$.

Suppose that $h_i \in \ideal(x_1,x_2,x_3^2)$ for each $i \ge 3$.
Then $h_i^{*} \in \ideal(x_1,x_2,x_3)$ for each $i \ge 3$. If we combine
this with $h_i^{*} \in K$ for every $i \in \{3,4,\ldots,k\}$, we infer
that $h_i^{*} = 0 = (h_i)_{x_3}$ for all $i \in \{3,4,\ldots,k\}$.
Hence $(H_1)_{x_1} + (H_2)_{x_2} = 0$ follows from equation
\eqref{eq2.4a} or equation \eqref{eq2.4b}. \qedhere

\end{enumerate}
\end{proof}

We conclude the proof of equation \eqref{eq2.4b} and Lemma \ref{lmeq2.3} (ii)
by describing an algorithm to find a $w$ and a linear transformation $T$ on the
third and subsequent coordinates, such that after replacing $H_1$ by $H_1(Tx)$,
the condition of Lemma \ref{histar} is met. The values of $w(x_i)$ will be
real (rational) numbers.

Initially, we take $w(x_3) = w(x_4) = \cdots = w(x_n) = 1$. This yields
a polynomial $\Phi(H_1)$, and by way of linear transformation on the third
and subsequent coordinates, we can obtain that for some $k \ge 3$,
\begin{gather*}
\Phi(H_1)_{x_3}, \Phi(H_1)_{x_4}, \ldots, \Phi(H_1)_{x_k}
~\mbox{are linearly independent over $K$}, \\
\Phi(H_1)_{x_{k+1}} = \Phi(H_1)_{x_{k+2}} = \cdots = \Phi(H_1)_{x_n} = 0, \\
1 = w(x_3) \le w(x_4) \le \cdots \le w(x_k) = w(x_{k+1}) = w(x_{k+2}) = \cdots = w(x_n).
\end{gather*}
If $H_1 \in K[x_1,x_2,x_3,\ldots,x_k]$, then we are done, so assume otherwise.
We show that it is possible to increase $k$.

We replace $w(x_{k+1}) = w(x_{k+2}) = \cdots = w(x_n)$ by a larger real number
$\gamma$, in such a way that $\Phi(H_1)$ gets other terms besides those it already had.
Suppose that
$$
c_3 \Phi(H_1)_{x_3} + c_4 \Phi(H_1)_{x_4} + \cdots + c_n \Phi(H_1)_{x_n} = 0
$$
for certain $c_i \in K$, not all zero. By splitting the left-hand side
in a part of monomials $\tau$ for which $w(\tau) \le w(H_1) - \gamma$
and another part of monomials $\tau$ for which $w(\tau) > w(H_1) - \gamma$,
we infer that $c_3 = c_4 = \cdots = c_k = 0$, and
$$
c_{k+1} \Phi(H_1)_{x_{k+1}} + c_{k+2} \Phi(H_1)_{x_{k+2}} + \cdots + c_n \Phi(H_1)_{x_n} = 0.
$$
Just as in the initial case, we can transform all linear dependencies
between $\Phi(H_1)_{x_{k+1}}, \Phi(H_1)_{x_{k+2}}, \ldots, \Phi(H_1)_{x_n}$
to
$$
\Phi(H_1)_{x_{k'+1}} = \Phi(H_1)_{x_{k'+2}} = \cdots = \Phi(H_1)_{x_n} = 0 \\
$$
for some $k' \ge k$.
But $k' > k$ because $\Phi(H_1) \notin K[x_1,x_2,x_3,\ldots,x_k]$.
So we can increase $k$ until $H_1 \in K[x_1,x_2,\ldots,x_k]$, which yields
the condition of Lemma \ref{histar}.

\end{document}